\documentclass[12pt,leqno]{article} 
\usepackage{amsmath,amsfonts,amscd,Martintxo,url,txfonts}
\usepackage{graphicx}
\usepackage{color}
\usepackage[all]{xy}

\newcommand{\refbat}{\arabic{Refbat}}  
\title{\Large\bf Finitness of the basic intersection cohomology of a Killing foliation}

\author{Martintxo Saralegi-Aranguren\thanks{Univ Lille Nord de France F-59\kern 1mm 000 LILLE, FRANCE.
UArtois, Laboratoire de MathŽmatiques de Lens EA~2462.
FŽdŽration CNRS Nord-Pas-de-Calais FR~2956.
Facult\'e des Sciences Jean Perrin.
Rue Jean Souvraz, S.P.\kern 1mm 18.
 F-62\kern 1mm 300 LENS, FRANCE.
{\sl saralegi@euler.univ-artois.fr}. }
\\ {\small UniversitŽ d'Artois } 
\and   
Robert Wolak\thanks{Instytut Matematyki. Uniwersytet Jagiellonski. 
Stanislawa Lojasiewicza 6, 30\kern 1mm 348
KRAKOW, POLAND.
{\sl robert.wolak@im.uj.edu.pl.}  Partially supported by the KBN 
grant 2 PO3A  
021 25.}
\\ {\small Uniwersytet Jagiellonski}
}

\begin{document}  
\maketitle

\begin{abstract}  We prove that the basic intersection cohomology 
$\lau{\IH}{*}{\per{p}}{\mf}$, where $\F$ is  the singular 
foliation  determined by an isometric action 
of a 
Lie group $G$ on the compact manifold $M$, is finite dimensional.
 \end{abstract}

\bigskip

This paper deals with an action $\Phi \colon G \times M \to M$ of a  
Lie group on a compact manifold preserving a riemannian metric on it. The 
orbits of this action define a singular foliation $\mathcal{F}$ on $M$. 
Putting together the orbits of the same dimension we get a stratification 
 of $M$. This structure is still very regular. The foliation 
$\mathcal{F}$ is in fact a 
conical foliation and we can define the 
basic intersection cohomology $\lau{\IH}{*}{\per{p}}{\mf}$ (cf. \cite{SW2}).
This invariant becomes the basic cohomology $\hiru{H}{*}{\mf}$ when the 
action $\Phi$ is almost free, and the intersection cohomology 
$\lau{\IH}{*}{\per{p}}{M/G}$ when the Lie group $G$ is compact.

The aim of this work is to prove that this cohomology  
$\lau{\IH}{*}{\per{p}}{\mf}$ is finite dimensional. This result generalizes \cite{EH} (almost free 
case), \cite{SW3} (abelian case) and \cite{SW2} (compact case).

\medskip

The paper is organized as follows. In Section 1 we present the 
foliation $\F$. The basic intersection 
cohomology $\lau{\IH}{*}{\per{p}}{\mf}$  associated to this foliation is studied in Section two. Twisted products are studied in Section 3.
The finiteness of $\lau{\IH}{*}{\per{p}}{\mf}$ is  proved in Section 4.

\medskip

In the sequel $M$  is a connected, second countable, Haussdorff,   
without boundary and smooth 
(of class $C^\infty$) manifold  of dimension
$m$. 
All the maps are considered smooth unless something else is indicated.

\section{Killing foliations determined by isometric actions.}
\label{KF}
We study in this work the foliations induced by isometric actions: the {\em Killing foliations}. These foliations are examples of the conical foliations for which the  basic intersection cohomology  has been defined (see \cite{SW2,SW3}).
We present this geometrical framework in this section.

 \smallskip
 
 \prg {\bf Killing foliations}. A smooth action $\Phi \colon G \times M 
 \to M$ of a Lie group $G$ on a manifold $M$ is a {\em isometric action} 
when there exists a riemannian metric $\mu$ on $M$ preserved by $G$.
 
  The connected components of the orbits of the action $\Phi$  determine 
  a partition $\mathcal{F}$ 
  on $M$. In fact, this partition is a singular riemannian foliation 
  that we shall call {\em Killing foliation} (cf. \cite{Mo}). Notice that $\F$ is also a conical foliation in the sense of \cite{SW2,SW3}.
 Classifying the points of $M$ following the dimension of the leaves of $\F$
one gets the {\em stratification} $\SF$ of $\mathcal{F}$. It is 
determined by the equivalence relation
$
x \sim y \Leftrightarrow \dim G_{x} = \dim G_{y}.$
The elements  of  $\SF$ are called {\em strata}.

In the particular case where  the closure of $G$ in the
 isometry group of 
 $(M,\mu) $ is a compact Lie group\footnote{This is always the case when the manifold 
 $M$ is a compact.}  we shall say that the action $\Phi$ is a {\em tame action}. In fact, a smooth action $\Phi \colon G \times M \to M$ is tame 
 if and only if it extends to a smooth action $\Phi \colon K \times M 
 \to M$ where $K$ is a compact Lie group containing $G$ (cf. \cite{K}). 
The group $K$ is not unique, but we always can choose $K$ in such a 
  way that $G$ is dense in $K$. We shall say that $K$ is a {\em tamer group}. 
Here the strata  of  $\SF$ are $K$-invariant closed submanifolds of $M$.

Since the aim of this work is the study of $\mathcal{F}$ and not the action $\Phi$ itself,  we can consider that the Lie group $G$ is
connected. Let us see that.
 \bP
 \label{normal}
 Let $\Phi \colon G \times M \to M$ is a tame action.
Let $G_0$ be  the connected component of $G$ containing the unity element. The Killing folation 
 defined by the restriction $\Phi \colon G_0 \times M \to M$ is also
 $\mathcal{F}$.
 \eP
  \pro
  The partition 
   $\mathcal{F}$ is defined by 
  this equivalence relation:
  $$
  x \sim y \Longleftrightarrow \exists \hbox{ continuous path } \alpha \colon 
  [0,1] \to G(x)  \hbox{ such that } \alpha(0) = x \hbox{ and }  \alpha(1) = 
  y.
  $$
  Since the map $\Delta  \colon G \to G(x)$, defined by $\Delta (g) = \Phi(g,x) = g \cdot x$,  is  a submersion (see for example \cite{Br}) then  
\bee
  x \sim y \Longleftrightarrow \exists \hbox{ continuous path } \beta \colon 
  [0,1] \to G  \hbox{ such that } \beta(0) = e \hbox{ and }  \beta(1) \cdot x= 
  y,
  \eee
   and by definition of $G_0$
  \bee
  x \sim y \Longleftrightarrow \exists \hbox{ continuous path } \beta \colon 
  [0,1] \to G_0  \hbox{ such that } \beta(0) = e \hbox{ and }  \beta(1) \cdot 
  x= 
  y.
\eee
  This gives the result. \qed

When $G$ is connected, the tamer group $K$ has richer properties.
 \bP
   \label{trio}
Let $G$ be a connected Lie subgroup of a compact Lie group $K$. If $G$ is dense in $K$ then $G \triangleleft K$ and the quotient group $K/G$ is commutative.
  \eP
\pro
The Lie algebra $\mathfrak{g}$ is $\Ad_G$- invariant and hence, by density, $\Ad_K$-invariant. Then $\mathfrak{g}$ is an ideal of $\mathfrak{k}$. The connectedness of $G$ gives that $G$ is a normal subgroup of $K$. Since $\Ad_G$ acts trivially on $\mathfrak{k}/\mathfrak{g}$, $\Ad_K$ acts trivially, too. Therefore, $\mathfrak{k}/\mathfrak{g}$ is abelian (see for example \cite[pag. 628]{P}). \qed

 \prg 
  \addtocounter{Sei}{\value{proposition}} 
{\bf Particular  tame actions.} A  {\em trio} is a triple $(K,G,H)$, with $K$ is a compact Lie group, $G$ a normal subgroup of $ K$ and $H$  a closed subgroup of $K$. 
We present now some tame actions associated to a trio  $(K,G,H)$. They are going to be intensively used in this work. First of all we need some definitions.

 \begin{enumerate}
\item[-]  The action $\Phi_l \colon K \times K \to K$ is defined by  $\Phi_l(g,k) = g \cdot 
 k $. For each element $u$ of the Lie algebra  $\mathfrak{k}$ of 
$K$, we shall write $X^u$ the associated (right invariant) vector field. It is 
defined by $X^{u}(k) = T_{e}R_{k}(u)$ where $R_{k} \colon K 
\to K$ is given by $R_{k}(\ell) = \ell \cdot k$. 

\item[-]  The action $\Phi_r \colon K \times K \to K$ is defined by  
$\Phi_r(g,k) = k\cdot 
g^{-1}$. For each element $u \in \mathfrak{k}$ of 
$K$, we shall write $X_u$ the associated (left invariant) vector field. It is 
defined by $X_{u}(k) = -T_{e}L_{k}(u)$ where $L_{k} \colon K 
\to K$ is given by $L_{k}(\ell) = k \cdot \ell $.

\item[-]  The action $\Psi \colon  K \times K/H \to K/H$ is defined 
 by  $\Psi(g,kH) = (g\cdot 
 k)H$. For each element $u \in\mathfrak{k}$, we shall write $Y_u $
 the associated vector field. Since the canonical projection $\pi \colon K \to K/H$ is a $K$-equivariant map, 
then we have $\pi_{*}X^{u} = 
Y_{u}$ for each $u \in \mathfrak{k}$.

 \item[-] The action $\Gamma \colon  H \times H \to H$ is defined 
 by  $\Gamma(g,h) = g \cdot 
 h$. For each element $u$ of the Lie algebra $\mathfrak{h}$
of $H$  we write $Z^u$ the associated (right invariant) vector field.

\end{enumerate}

\nt The associated actions we are going to use are the following.

 \Zati {\em The restriction $\Phi_l \colon G \times K \to K$, which  induces the regular Killing  foliation  $ \mathcal{K} $}. 

 \medskip

\zati {\em The restriction $\Phi_r \colon G \times K \to K$, which  induces the regular Killing  foliation  $ \mathcal{K} $}. 

\medskip

\nt  
Since $G \triangleleft K$, the foliation $ \mathcal{K} $ is determined by the family of vector fields 
$\{ X^u \ / \ u \in \mathfrak{g}\}$, where $\mathfrak{g}$ is the Lie 
algebra of $G$, and also by the family 
$\{ X_u \ / \ u \in \mathfrak{g}\}$.
The orbits $G(k) = 
 G   k = k  G$ 
 have the same dimension $\dim G$.

 \zati {\em 
 The restriction $\Psi \colon G \times K/H \to K/H$,  which induces the regular Killing foliation $ \mathcal{D} $}. 

 \medskip

 \nt 
The foliation $ \mathcal{D} $ is determined by the family of vector fields 
$\{ Y_u \ / \ u \in \mathfrak{g}\}$.   The orbits $G (kH)  $ have the same dimension 
 $\dim G - \dim (G \cap H)$.

 \zati {\em The restriction $\Gamma \colon  (G \cap H) \times H \to H$, which induces the regular Killing foliation $ \mathcal{C} $}.

 \medskip 
 
 \nt  
 The foliation $ \mathcal{C} $ is determined by the family of vector fields 
 $\{ Z^u \ / \ u \in \mathfrak{g} \cap \mathfrak{h}  \}$. 
 The orbits $(G \cap H)(k) $ have the same dimension 
  $\dim( G \cap  H)$. 

 \medskip 

  \zati {\em The restriction $\Phi_r\colon GH  \times K \to K$,  which induces the regular Killing
  foliation $ \E $}.
 
  \medskip 
  
  \nt 
Notice that $GH$ is  a Lie group since $G$ is normal in $K$.
    The foliation $ \E $ is, in fact,  determined by the vector fields 
   $\{ X_u \ / \ u \in \mathfrak{g} + \mathfrak{h}  \}$.  The orbits $(G  H)(k) $ have the same dimension 
  $\dim G + \dim H -\dim( G \cap  H)$. 
   
 \prg {\bf Twisted product.} 
  \addtocounter{Due}{\value{proposition}} 
  In order to prove the finiteness of the basic intersection cohomology we decompose the manifold in a finite number of simpler pieces. These are the twisted products we introduce now. 
 
 We fix a trio $(K,G,H)$  and a smooth action $\Theta \colon H \times N \to N  $ of $H$ on the manifold $N$.   The {\em twisted 
   product} is the quotient
   $K \bi{\times}{H}N $  of $K \times N $ by the 
   equivalence relation
   $ (k,z) \sim (k\cdot h^{-1},\Theta(h,z) = h \cdot z) .
   $
   The element of   $K \bi{\times}{H}N $  corresponding 
   to $(k,z) \in K \times N $ is denoted by 
   $ <k,z>$. This
   manifold is endowed with the tame action
   $$
   \Phi \colon G \times (K \bi{\times}{H}N ) \TO (K \bi{\times}{H}N ),
   $$
   defined by $\Phi (g, <k,z>) = <g \cdot k,z>$. We denote by
   $\mathcal{W}$ the induced Killing
  foliation.

    We also use the following tame action, namely, 
   the restriction 
   $$\Theta \colon (G \cap H )\times N \to N $$ 
   whose 
   induced Killing foliation is denoted by $\mathcal{N}$.

   The canonical projection $\Pi \colon K \times N \to K \times_{H} N $ relates the involved foliations as follows:
   
   \Zati 
   \addtocounter{Sette}{\value{zatia}} 
   $\Pi_* (\mathcal{K}\times \mathcal{I}) = \mathcal{W}$, where $\mathcal{I}$ is the pointwise foliation (since the map $\Pi$ is $G$-equivariant).

   \zati 
    \addtocounter{Tre}{\value{zatia}} 
   $\Str{\mathcal{W}}= 
     \{ \Pi(K   \times S) \ / \  
    S \in \SN\} = \Pi \left(\{ K \}Ê\times \SN
    \right)$ (since $G_{<k,z>} = k(G \cap H)_{z}k^{-1}$).

\section{ Basic Intersection cohomology}
\label{bs}
In this section we recall the definition of the basic intersection\footnote{We refer the reader to   \cite{SW2},\cite{SW3} for details.}  
cohomology and we present the main properties we are going to use in this work.
For the rest of this section, we fix a conical foliation $\mathcal{F}$ defied on a manifold $M$. The associated stratification is $\SF$. The regular stratum of is denoted by $\bi{R}{\F}$. We shall write $m = \dim M$, $r = \dim \mathcal{F}$ and $ s=m - r = 
\codim_{M}\mathcal{F}$.

\smallskip

We are going to deal with 
differential forms on  a product $(\hbox{manifold}) \times [0,1[^p$, they 
are restrictions of differential forms defined on $(\hbox{manifold})\times 
]-1,1[^p$.

\prg{\bf Perverse forms.} 
  \addtocounter{Nove}{\value{proposition}} 
Recall that a {\em conical chart} is a foliated diffeomorphism $\phii \colon  (\R^{m-n-1}  \times c \S^n, \mathcal{H} \times c \G) \to  (U,\F_U) $ where $ (\R^{m-n-1} , \mathcal{H})$ is a simple foliation and $(\S^n,\G)$ is a conical foliation without  $0$-dimensional leaves. We also shall denote this chart by $(U,\varphi, S)$ where $S$ is the stratum of $\SF$ verifying  
$\phii( \R^{m-n-1}  \times \{ \vartheta \}) = U \cap S$.

The differential complex  $\lau{\Pi}{*}{\mathcal{F}}{M \times 
[0,1[^p}$ 
of {\em perverse forms} of $M \times 
[0,1[^p$ is introduced by induction 
on  $\depth  \SF$. 
When this depth is  0 then 
$$\lau{\Pi}{*}{\mathcal{F}}{M \times [0,1[^p} = 
\hiru{\Omega}{*}{M 
\times [0,1[^p}.
$$

Consider now the 
generic case. 
A perverse form of $M \times 
[0,1[^p$ is  first of all a differential form $\om\in \hiru{\Om}{*}{R_\mathcal{F} \times [0,1[^p}$ such that,  
\begin{equation*}
\left\{
\begin{array}{ll}
    \hbox{the pull-back \hspace{2cm}} &
 (\phii \times \mathbb{I}_{[0,1[^p})^*\omega \in \hiru{\Om}{*}{\R^{m-n-1} \times \bi{R}{\G}\times ]0,1[ \times [0,1[^p}\\[.3cm]
\hbox{extends to} &
\bi{\om}{\varphi} \in \lau{\Pi}{*}{\mathcal{H}\times c\G}{\R^{m-n-1} \times \S^n \times [0,1[^{p+1}}
\end{array}
\right.
\end{equation*}
for  any conical chart  $(U,\phii)$, where $\bi{\mathbb{I}}{\bullet}$ stands for the identity map.
Notice that $\hiru{\Om}{*}{M}$ is included on $\lau{\Pi}{*}{\mathcal{F}}{M}$\footnote{Through the restriction $\om \mapsto 
\bi{\om}{\bi{R}{\F}}$.}.

\prg{\bf Perverse degree.} 
 \addtocounter{Dieci}{\value{proposition}} 
The amount of transversality of a perverse form $\om \in \lau{\Pi}{*}{\mathcal{F}}{M}$ 
with respect to a  
singular stratum $S \in \SF$
is measured by the perverse degree $\bi{||\om||}{S}$. We recall here the definition of local 
perverse degree $\bi{||\om||}{U}\in  \{ -\infty\} \cup \N  $ of $\om$ relatively to a conical chart $(U,\phii,S)$:

\begin{itemize}

	\item[1.] $\bi{||\om||}{U}= -\infty$ when $\bi{\om}{\varphi} \equiv 0$ on 
		   $\R^{m-n-1} \times \bi{R}{\G}  \times \{ 0 \}$,
	\item[2.]  $\bi{||\om||}{U}\leq p$, with  $p \in \N$,  when $\bi{\om}{\varphi} 
	(v_{0}, \ldots ,v_{p},-)\equiv 0$ where the vectors $\{ v_{0}, 
	\ldots , v_{p}\}$ are tangent to the fibers of 
	$\bi{P}{\varphi} \colon \R^{m-n-1} \times \bi{R}{\G}  \times \{ 0 \}  \TO U \cap S$\footnote{The map $\bi{P}{\varphi} \colon \R^{m-n-1} \times \S^{n} \times [0,1[  \TO U$ is defined by $\bi{P}{\varphi} (x,y,t) = \phii (x,[y,t])$.}.

\end{itemize}
This number does not depend on the choice of the conical 
chart  (cf. \cite[Proposition 1.3.1]{SW3}). 
Finally, we define the {\em 
perverse degree} 
$\bi{||\omega||}{S}$ by
\begin{equation*}
\bi{||\omega||}{S} =  \sup \left\{ \bi{|| \om ||}{U} \ / \  (U,\phii,S) \hbox{ 
conical chart }  \right\}.
\end{equation*}
  The perverse degree of  $\om \in \hiru{\Om}{*}{M}$ verifies $\bi{||\omega||}{S} \leq 0$ for any singular stratum $S \in \SF$ (cf. \ref{bs}.\nove). 

\prg {\bf Basic cohomology}. 
\setcounter{Refbat}{\value{proposition}}  The basic cohomology of the foliation 
$\mathcal{F}$  is an important  tool 
to study its transversal structure and plays the r™le of the cohomology 
of the 
orbit space $\mf$, which can be a wild topological espace. 
A differential form $\om \in \hiru{\Om}{*}{M}$ is {\em basic} if 
$
i_X \om = i_X d \om = 0,
$
for each vector field $X$ on $M$ tangent to the foliation $\mathcal F$. For exemple, a 
function $f$ is basic iff $f$ is constant on the leaves of $\F$. We shall write 
$\hiru{\Om}{*}{\mf}$ for
the complex of basic forms.  Its cohomology $\hiru{H}{*}{\mf}$ is the {\em 
basic cohomology} 
of $(M,\mathcal{F})$.
We also use  the {\em relative basic cohomology} $\hiru{H}{*}{(M,F)/
\mathcal{F}}$, that is, the cohomology computed from the complex of basic forms vanishing 
on the saturated set $F \subset M$.
 The basic cohomology does not use the stratification $\SF$. 

\prg {\bf Basic intersection cohomology}.  
 A {\em perversity} is a map $\per{p} 
\colon \SiF \to \Z \cup \{-\infty,\infty\}$, where $\SiF$ is the family of singular strata. 
The {\em constant perversity} $\per{\iota }$ is defined by $\per{\iota }(S) 
=\iota $, where $\iota \in \Z \cup \{-\infty,\infty\}$.

The basic intersection cohomology appears when one considers basic perverse forms whose 
perverse degree is controlled by a perversity.
We shall put
$$
\lau{\Om}{*}{\per{p}}{\mf}  = \left\{ \om \in 
\lau{\Pi}{*}{\mathcal{F}}{M} \ / \  \om \hbox{ is basic and }
\max \left(\bi{||\om||}{S} ,\bi{||d\om||}{S}\right)\leq \per{p}(S) 
\ \ \forall S \in \SiF \right\}
$$
the complex of basic perverse forms whose perverse degree (and that of the their 
derivative) is bounded by the perversity $\per{p}$. The cohomology $\lau{\IH}{*}{\per{p}}{\mf}$ of this complex is the {\em 
basic intersection cohomology\footnote{BIC for short.}} of $(M,\mathcal{F})$ relatively to the perversity $\per{p}$.

   Consider a twisted product   $K \bi{\times}{H}N $.
Perversities on $K \bi{\times}{H}N $ and $ K \times N $ are determinate by perversities on $N $ by the formula  (cf. \ref{KF}.\due \tre):
\be
\label{pe}
\per{p}(K \times S) =  \per{p}(\Pi(K \times S)) = 
\per{p}(S).
\ee

\prg {\bf Mayer-Vietoris}. 
  \addtocounter{Quince}{\value{proposition}} 
This is the technique we use in order to decompose the manifold in nicer pieces. An open covering $\{ U , V \}$  of   $M$  by saturated open 
subsets is a {\em basic covering}. It possesses a 
subordinated partition of the unity made up of  basic  functions defined 
on $M$ (see \cite{RSW1}).
For a such covering we have the Mayer-Vietoris short sequence
$$
0 \to \lau{\Om}{*}{\per{p}}{\mf} \to 
\lau{\Om}{*}{\per{p}}{U/\mathcal{F}} \oplus 
\lau{\Om}{*}{\per{p}}{V/\mathcal{F}}  \to
\lau{\Om}{*}{\per{p}}{(	U \cap V)/\mathcal{F}}  \to 0,
$$
where the map are defined by $\om \mapsto (\omega,\omega)$ and 
$(\alpha,\beta) \mapsto \alpha - \beta$. The third map is onto since the 
elements of the partition of the unity are {\em controlled functions}, id est, elements of $\lau{\Om}{0}{\per{0}}{-}$ (cf. \ref{bs}.\dieci). Thus,  the sequence 
 is exact. This result is not 
longer true for 
more general coverings.

We shall use in this work the two following local calculations (see \cite[Proposition 3.5.1 and Proposition 3.5.2]{SW3} 
for the proofs).

\bp
\label{pro}
Let $(\R^k,\mathcal{H})$ be a simple foliation. Consider  $\per{p}$ a 
perversity on $M$ and define the perversity $\per{p}$ on $\R^k \times 
M$  by $\per{p}(\R^k \times S) = \per{p}(S)$.
The canonical projection $\pr \colon \R^k \times M \to M$ induces the isomorphism
$$
\lau{\IH}{*}{\per{p}}{\mf} \cong \lau{\IH}{*}{\per{p}}{\R^k \times M /  \mathcal{H} 
\times \mathcal{F}}.
$$
\ep

\bp
\label{local5}
Let $\G$ be a conical foliation without 0-dimensional leaves on 
the sphere $\S^n$. A perversity 
$\per{p}$  on $c\S^n$ gives the perversity $\per{p}$ on $\S^n$ defined 
by $\per{p}(S) = \per{p}(S \times ]0,1[)$. The canonical projection 
$\pr\colon 
\S^n\times 
]0,1[ \to \S^n $ induces the isomorphism
$$
\lau{\IH}{i}{\per{p}}{c\S^n / c\G }=
\left\{
\begin{array}{cl}
\lau{\IH}{i}{\per{p}}{\S^n / \G }& \hbox{if } i \leq \per{p}(\{ 
\vartheta \}) \\
0 & \hbox{if } i > \per{p}(\{ \vartheta \}). 
\end{array}
\right.
$$
\ep

In the next section we shall need the following  technical Lemma.

\bl
\label{iX}
Let $\Phi \colon K \times M \to M$ be a smooth action, where $K$ is a compact Lie group, and let $V$ be a
fundamental vector field of this  action. Consider  a  normal subgroup $G$ of $K$ and write $\F$ the associated conical foliation on $M$.
Then, the interior  operator
$
i_{V} \colon \lau{\Om}{*}{\per{p}}{\mf} \TO \lau{\Om}{*-1}{\per{p}}{\mf}
$
is well defined, for any perversity $\per{p}$.
\el

\pro 
Since the question is a local one, then it suffices to consider where $M$ is a twisted product
$K \bi{\times}{H}N$\footnote{In fact, $N$ is an euclidean space $\R^a$ et $\Theta$ is an orthogonal action.}. 
Notice that the blow up $\Pi 
\colon K \times N \to K \bi{\times}{H}N $ is a $K$-equivariant map relatively to the action
$\ell \cdot (k,z) = (\ell \cdot k,z)$. This gives $\Pi_{*}(X^{u},0) = V$ for some $u \in \mathfrak{k}$.
From Lemma \ref{twist} we know that it suffices to prove that the operator
$$
i_{ (X^u,0) } \colon \lau{\Om}{*}{\per{p}}{K \times N / \mathcal{K}  \times  \mathcal{N}} \TO
\lau{\Om}{*-1}{\per{p}}{K \times N / \mathcal{K}  \times  \mathcal{N}}
$$
is well defined.  Since $G \triangleleft K$ then the vector field $X^u$ preserves the foliation $\mathcal{K}$. So, it suffices to prove that 
 the operator
$$
i_{ (X^u,0) } \colon \lau{\Om}{*}{\per{p}}{K \times N } \TO
\lau{\Om}{*-1}{\per{p}}{K \times N }
$$
is well defined.  This comes from the fact that $X^u$ acts on the $K$-factor while the perversion conditions are measured on the $N$-factor
(cf. \eqref{pe}).\qed

\section{The BIC of a twisted product}
\label{tp}
 We compute now the BIC of a twisted product $K \bi{\times}{H}N $ (cf. \ref{KF}.\due ) for a perversity $\per{p}$ (cf. \eqref{pe}).

\bl
\label{twist}
The natural projection 
$\Pi \colon K \times N  \to K  \bi{\times}{H} N $ induces the differential 
monomorphism
\begin{equation}
\label{w}
\Pi^* \colon \lau{\Om}{*}{\per{p}}{K \bi{\times}{H}N / \mathcal{W}}
\TO 
\lau{\Om}{*}{\per{p}}{K \times N /\mathcal{K}\times 
 \mathcal{N}}.
\end{equation}
 Moreover, given a differential form $\om$ on 
$K \bi{\times}{H} R_{\mathcal{W}}$, we have:
\begin{equation}
\label{ww}
\Pi^* \om \in \lau{\Om}{*}{\per{p}}{K \times N /\mathcal{K}\times
 \mathcal{N}} \Longleftrightarrow \om \in \lau{\Om}{*}{\per{p}}{K \bi{\times}{H}N / \mathcal{W}}.
\end{equation}
\el
\pro
    Notice that the injectivity of $\Pi^{*}$ comes from the fact that $\Pi$ is a 
   surjection. For the rest, we proceed in  several steps.
   
   \Zati \underline{\em A  foliated atlas for 
    $\pi \colon K \to 
   K/H$}.

   Since $\pi \colon K \to K/H$ is a $H$-principal bundle then it possesses an atlas 
   $\mathcal{A} =\left\{ 
   \phii \colon \pi^{-1}(U) \TO U \times H \right\}$ made up with $H$-equivariant 
   charts: $\phii (k \cdot h^{-1}) = (\pi (k), 
   h\cdot h_0)$ if $\phii (k ) = (\pi (k), 
   h_0)$. We study the foliation $\phii_*\mathcal{K}$.
   This equivariance property gives $\phii_* X_u = (0,Z^u)$ for each $u 
   \in \mathfrak{g} \cap \mathfrak{h}$. 
   Thus, the trace of the foliation $\phii_*\mathcal{K}$ on the 
   fibers of the canonical projection $\pr\colon U \times H \to U$ is 
   $\mathcal{C}$.
   On the other hand, since the map $\pi$ is a $G$-equivariant map then $\pi_* 
   \mathcal{K}=  \mathcal{D} $, which gives $\pr_*\phii_*\mathcal{K}= 
    \mathcal{D} $. We conclude that $\phii_*\mathcal{K}\subset  \mathcal{D}  \times 
   \mathcal{C}$. By dimension reasons we get $ \phii_*\mathcal{K}=  \mathcal{D}  \times 
   \mathcal{C}$.
   The atlas $\mathcal{A}$ is an $H$-equivariant foliated atlas of 
   $\pi$.

   \zati \underline{\em A foliated atlas for $ \Pi \colon K \times N \to K 
\bi{\times}{H}N $}.

   We claim that  $\mathcal{A}_{_{\#}} = 
   \left\{\overline{\phii} \colon \pi^{-1}(U) \bi{\times}{H}N  \TO U \times 
   N 
   \ / \ (U,\phii) \in \mathcal{A} \right\}$ is a foliated atlas of 
   $K 
\bi{\times}{H}N $
   where
   the map 
   $
   \overline{\phii} 
   $
   is defined by $\overline{\phii}(<k,z>) =  
   (\pi(k),(\Theta ( (\phii^{-1}(\pi(k),e))^{-1}  \cdot k ,  z)))$. This 
   map  is a 
   diffeomorphism whose inverse 
   is $\overline{\phii}^{-1}(u,z) = <\phii^{-1}(u,e),z>$.  It verifies
   $$
    \overline{\phii}_*\mathcal{W} \stackrel{\ref{KF}. \due \sette}{= \! =\! =}  \overline{\phii}_*\Pi_* ( \mathcal{K}  
   \times \mathcal{I}) = 
   \overline{\phii}_*\Pi_* (\phii^{-1} \times \Ide_{N })_*(
    \mathcal{D}  \times 
\mathcal{C}\times \mathcal{I}).
   $$
   A straightforward calculation shows
   $\overline{\phii} \rondp \Pi \rondp (\phii^{-1} \times \Ide_{N }) = 
   (\Ide_U \times \Theta)$. 
   Since $\mathcal{C}$ is defined by the action $\Gamma$ then $\Theta_* (\mathcal{C} \times \mathcal{I}) = 
    \mathcal{N}$.
   Finally we obtain
   $
   \overline{\phii}_*\mathcal{W}=   \mathcal{D}  \times  \mathcal{N}.
   $

   \zati \underline{\em Last Step}.

    Given $(U,\phii) \in \mathcal{A}_{\#}$, we have the commutative diagram
    \begin{equation*}
\xymatrix@C=2.5cm{
 \ar[d]_{Q}  U \times H \times N  \ar[r]^{\varphi^{-1} \times \Ide_{N }} & 
\ar[d]^\Pi K  \times N  \\
   U \times N  \ar[r]^{\overline{\varphi}^{-1} } &K \bi{\times}{H}N 
}
\end{equation*}
   where $Q(u,h,z) = (u, h^{-1} \cdot z)$,
   $\Pi^{-1}(\Ima \overline{\phii}^{-1}) = \Ima \left(\phii^{-1} \times \Ide_{N }\right)$ and the rows are foliated imbeddings.
  Now, since \eqref{w} and \eqref{ww} are local questions  
   then it suffices to 
   prove that 
    \begin{itemize}
\item[-]  $Q^* \colon \lau{\Om}{*}{\per{p}}{U \times N  /  \mathcal{D}  \times 
    \mathcal{N}}
   \TO 
   \lau{\Om}{*}{\per{p}}{U \times H \times N /  \mathcal{D}  \times 
\mathcal{C}\times 
    \mathcal{N}}
   $
   is well-defined, and 
   \item[-] $Q^* \om \in  \lau{\Om}{*}{\per{p}}{U \times H \times N /  \mathcal{D}  \times 
\mathcal{C}\times 
    \mathcal{N}}\Longleftrightarrow \om \in \lau{\Om}{*}{\per{p}}{U \times N  /  \mathcal{D}  \times 
    \mathcal{N}}$, for any $\om \in  \hiru{\Om}{*}{U \times R_{ \mathcal{N}}}$.
\end{itemize}
  This comes from the fact that the map 
   $$\nabla \colon (U \times H \times N ,  \mathcal{D}  \times 
\mathcal{C}\times 
    \mathcal{N}) \TO (U \times H \times N ,  \mathcal{D}  \times 
\mathcal{C}\times 
    \mathcal{N}),
   $$
   defined by $\nabla(u,h,z) = (u,h,h^{-1} \cdot z))$, is a 
   foliated diffeomorphism\footnote{Since $G \cap H \triangleleft H$.} and $Q =   \pr_0\rondp \nabla$, with $\pr_0 \colon 
   U \times H \times N  \to U \times N $ canonical projection
   (cf. Proposition \ref{pro}).
    \qed

\prg{\bf The Lie algebra  \mbox{\boldmath$ \frak{k}$}.}
        \addtocounter{Otto}{\value{proposition}} We suppose in this paragraph that that $G$ is also dense on $K$.
Choose $\nu$ a bi-invariant riemannian 
metric on $K$, which exists by compactness.
Consider
\begin{equation*}
B = \left\{ u_1, \ldots u_a,u_{a+1}, \ldots , u_{b}, u_{b+1},\ldots , 
u_c,u_{c+1},\ldots,u_f\right\}
\end{equation*}
an orthonormal basis of the Lie algebra $\mathfrak{k}$ of $K$ 
with $\left\{ u_1, \ldots u_b\right\}$ basis of the Lie algebra $\mathfrak{g}$ of $G$ and 
$\left\{ u_{a+1}, \ldots u_c\right\}$ basis of the Lie algebra $\mathfrak{h}$ of $H$. 
For each indice $1 \leq i \leq f$ 
we shall write $X_{i}\equiv X_{u_i}$  and $X^{i}\equiv X^{u_i}$ (cf. \ref{KF}.\hspace{-.1cm} \sei).

Let  $\gamma_i \in \hiru{\Om}{1}{K}$ be the 
dual form of $X_i$, that is, $\gamma_i = i_{X_i} \nu$.  Notice that $\delta_{ij} = 
\gamma_{j}(X_{i}) $.
These forms are invariant by the left
action of $K$.
Since the flow of $X^{j}$ is the multiplication on the left by $\exp (tu_j)$ then $L_{X^{j}} \gamma_{i} = 0$ for each $1 \leq j\leq f$.

For the differential, we have the formula
$
{\displaystyle d \gamma_l= \sum_{1\leq i < j \leq f} C^l_{ij} \gamma_i \wedge 
\gamma_j},
$
where 
$
{\displaystyle [X_i,X_j] = \sum_{l=1}^f C^l_{ij} X_l},
$
and $1 \leq i,j ,l\leq f$.
We have several restrictions on these coefficients.
Since $G \triangleleft K$ then   $\mathfrak{g}$ is an ideal of 
$\mathfrak{k}$ and therefore  we have 
\begin{equation*}
C^l_{ij}  = 0 \hbox{ 
for  $i \leq b < l$}. 
\end{equation*}
Since $K/G$ is an abelian group (cf. Proposition \ref{trio})  then the induced bracket on 
$\mathfrak{k}/\mathfrak{g}$ is zero and therefore we have
\begin{equation*}
C^l_{ij}  = 0 
\hbox{ 
for  } b < i,j,l \leq f.
\end{equation*}
These equations  imply that
\begin{equation}
\label{ciclo}
d\gamma_{l} = 0 
\hbox{ for each }b < l.
\end{equation}

The $ \E $-basic differential forms in $
\ib{\bigwedge}{*} 
(\gamma_{1}, \ldots,
\gamma_f)$ are exactly $
\ib{\bigwedge}{*}  
(\gamma_{c+1}, \ldots,\gamma_f)$ since they are cycles and 
the family $\{X_1,\ldots,X_c\}$ generates the foliation 
$ \E  $.
This gives
\begin{equation}
\label{KE}
   \hiru{H}{*}{K\Big /\mathcal{E}} = \ib{\bigwedge}{*} (\gamma_{c+1},\ldots,\gamma_f). 
\end{equation}

\prg
 \addtocounter{Doce}{\value{proposition}} 
{\bf Two actions of  \mbox{\boldmath$ H/H_0$}.} The Lie group $H$ preserves the foliation $ \mathcal{N}$ since the Lie group $G \cap H$ is a normal subgroup of $H$. Put $H_{0}$ the connected component of $H$ containing the unity element. Since it is a connected compact Lie group then   a standard argument shows that
\be
\label{hache}
\left(\lau{\IH}{*}{\per{p}}{ N /  \mathcal{N} } \right)^{H_{0}}=
\hiru{H}{*}{\left( \lau{\Om}{.}{\per{p}}{N /  \mathcal{N}}\right)^{H_{0}} }=
\lau{\IH}{*}{\per{p}}{N /  \mathcal{N} }
\ee
(cf. \cite[Theorem I, Ch. IV, vol. II]{GHV}). We conclude that the finite group $H/H_{0}$ acts naturally on $\lau{\IH}{*}{\per{p}}{N /  \mathcal{N} }$.

Since $H_0$ is a connected Lie subgroup of $GH$ then 
$
\left( \hiru{H}{*} {K / \mathcal{E} }  \right)^{H_{0}} = \hiru{H}{*} {K / \mathcal{E} } .
$
We conclude that the finite group $H/H_{0}$ acts naturally on $\hiru{H}{*} {K / \mathcal{E} } $.

\bp
Let $(K,G,H)$ be a trio with $G$ connected and dense in $K$. Then
\begin{equation*}
\lau{\IH}{*}{\per{p}}{K \bi{\times}{H} N / \mathcal{W} } = \left( \hiru{H}{*} {K / \mathcal{E} }  \otimes \lau{\IH}{*}{\per{p}}{N /  \mathcal{N}  }\right)^{H/H_{0}}.
\end{equation*}
\ep
\pro
Using the blow up $\Pi 
\colon K \times N \to K \bi{\times}{H}N $,
the computation of 
$\lau{\IH}{*}{\per{p}}{K \bi{\times}{H} N / \mathcal{W} } $ can be done by using the complex
$\Ima \left\{ 
\Pi^* \colon \lau{\Om}{*}{\per{p}}{K \bi{\times}{H}N / \mathcal{F}}
\TO 
\lau{\Om}{*}{\per{p}}{K \times N / \mathcal{K}  \times  \mathcal{N}}\right\}
$
(cf. Lemma \ref{twist}). We study this complex in several steps.
We fix
$
B = \left\{ u_1, \ldots ,u_f\right\}
$ an orthonormal basis of $\mathfrak{k}$ as in \ref{tp}.\otto.

\medskip

\nt $\langle  i \rangle $ \underline{Description of } $\hiru{\Om}{*}{K \times R_{ \mathcal{N}}}$. 

\smallskip

A differential form $\om \in 
\hiru{\Om}{*}{K \times R_{ \mathcal{N}}}$ 
is of the form
\begin{equation}
\label{form}
\eta + \sum_{1 \leq i_1 < \cdots < i_\ell\leq f}\gamma_{i_1} \wedge \dots 
\wedge \gamma_{i_\ell} \wedge \eta_{i_1, 
\ldots , i_\ell}, 
\end{equation}
where the forms $\eta, \eta_{i_1, 
\ldots , i_\ell}\in \hiru{\Om}{*}{K \times R_{ \mathcal{N}}}$ verify  
$i_{X_{j}} \eta = i_{X_{j}} \eta_{i_1, 
\ldots , i_\ell} = 0$ 
for each $1 \leq j \leq f$ and each 
$1 \leq i_1 < \cdots < i_\ell\leq f$.

\medskip

\nt $\langle  ii \rangle $  \underline{Description of } $\lau{\Pi}{*}{\mathcal{K}\times 
 \mathcal{N}}{K \times N }$. 

\smallskip

Since the foliation $\mathcal{K}$ is regular then we always can 
choose a conical chart of the form $(U_1 \times U_2, \phii_1 \times 
\phii_2)$ where $(U_1,\phii_1)$ is a foliated chart of $(K,
 \mathcal{K} )$ and $(U_2,\phii_2)$ is  a conical chart of 
$(N , \mathcal{N})$. The local blow up of the chart 
$(U_1 \times U_2, \phii_1 \times 
\phii_2)$ is constructed from the second factor without modifying the 
first one. 
So, the differential forms $\gamma_i$ are always 
perverse forms and a
differential form 
$\om \in \lau{\Pi}{*}{\mathcal{K}\times 
 \mathcal{N}}{K \times  N }$ 
is of the form \refp{form}
where $\eta, \eta_{i_1, \ldots , i_\ell}\in 
\lau{\Pi}{*}{\mathcal{K}\times 
 \mathcal{N}}{K \times  N }$ verify  
$i_{X_{j}} \eta = i_{X_{j}} \eta_{i_1, 
\ldots , i_\ell} = 0$ for each $1 \leq j \leq f$ and each 
$1 \leq i_1 < \cdots < i_\ell\leq f$.

\medskip

\nt $\langle  iii \rangle $ \underline{Description of } $\hiru{\Om}{*}{K \times R_{ \mathcal{N}}/
 \mathcal{K}  \times  \mathcal{N}}$. 

\smallskip

Take $\om \in 
\hiru{\Om}{*}{K \times R_{ \mathcal{N}}/ \mathcal{K}  \times 
 \mathcal{N}}$.
Since $\mathcal{K}$ is generated by the family $\{ X_{j} \ / \ 1 \leq j \leq b\}$ 
then $L_{X_{j}}\om = 0$ for any $1\leq j \leq b$, or equivalently, $R_{g}^{*}\om = \om$ 
for each $g \in 
G$ since $G$ is connected. By density, $R_{k}^{*}\om = \om$ for each $k \in 
K$ and therefore $L_{X_{j}}\om = 0$ for any $1 \leq j \leq f$ since $K$ is connected. 
We conclude that $L_{X_{j}}\eta 
= L_{X_{j}}\eta_{i_1, 
\ldots , i_\ell} = 0$ for any $1 \leq j \leq f$ and each 
$1 \leq i_1 < \cdots < i_\ell\leq f$.
This gives 
$
\om \in {\bigwedge}^*  
(\gamma_{1}, \ldots,\gamma_f)
\otimes \hiru{\Om}{*}{ R_{ \mathcal{N}}}.$

The $ \mathcal{N}$-basic differential forms of $
\hiru{\Om}{*}{ 
R_{ \mathcal{N}}}$ are exactly
$\hiru{\Om}{*}{ 
R_{ \mathcal{N}}/ \mathcal{N}}$.
The $\mathcal{K}$-basic differential forms of 
${\bigwedge}^*  
(\gamma_{1}, \ldots,\gamma_f)$ are exactly
${\bigwedge}^*  
(\gamma_{b+1}, \ldots,\gamma_f)$  (cf. \eqref{ciclo}).
From these two facts, we get
$$
\hiru{\Om}{*}{K \times R_{ \mathcal{N}}/
 \mathcal{K}  \times  \mathcal{N}}
=
{\bigwedge}^*  
(\gamma_{b+1}, \ldots,\gamma_f)
\otimes \hiru{\Om}{*}{ R_{ \mathcal{N}}/ \mathcal{N}}
$$
as differential graduate commutative algebras.

\medskip

\nt $\langle  iv \rangle $ \underline{Description of }
$\lau{\Om}{*}{\per{p}}{K \times N /  \mathcal{K}  \times 
 \mathcal{N}}$. 

\smallskip

From $\langle  ii \rangle  $ and $\langle  iii \rangle  $ it suffices to control the perverse degree of the 
forms 
$$
\eta + \sum_{b+1 \leq i_1 < \cdots < i_\ell\leq f}\gamma_{i_1} \wedge \dots 
\wedge \gamma_{i_\ell} \wedge \eta_{i_1, 
\ldots , i_\ell}  \in {\bigwedge}^*  
(\gamma_{b+1}, \ldots,\gamma_f)
\otimes \lau{\Pi}{*}{ \mathcal{N}}{ N }.
$$
 Consider $S$ a stratum of $\SN $. From $||\gamma_i||_{K \times S} = 0$ 
and  $||\eta||_{K \times 
S} =  ||\eta||_{S} $,  we get $||\gamma_{i_1} \wedge \dots 
\gamma_{i_\ell} \wedge \eta_{i_1, 
\ldots , i_\ell}||_{K \times S}  =||\eta_{i_1, 
\ldots , i_\ell}||_S$. We conclude that  
$$
\lau{\Om}{*}{\per{p}}{K \times N /  \mathcal{K}  \times 
 \mathcal{N}}
\cong
{\bigwedge}^*  
(\gamma_{b+1}, \ldots,\gamma_f)
\otimes 
\lau{\Om}{*}{\per{p}}{N /  \mathcal{N}}
$$
 (cf. \ref{KF}.\due \tre).  

\medskip

\nt $\langle  v \rangle $ \underline{Description of }
$\Ima \left\{ 
\Pi^* \colon \lau{\Om}{*}{\per{p}}{K \bi{\times}{H}N / \mathcal{F}}
\TO 
\lau{\Om}{*}{\per{p}}{K \times N / \mathcal{K}  \times  \mathcal{N}}\right\}
$.

\smallskip

We denote by
$\{W_{a+1},\ldots , W_c\}$ the
fundamental vector fields of the action $\Theta \colon H \times N \to 
N $ associated to the basis $\{u_{a+1}, \ldots , u_{c}\}$.
Consider now the action $\Upsilon\colon H \times (K \times N ) \to (K 
\times N )$ defined by $\Upsilon(h,(k,z)) = ( k \cdot h^{-1},\Theta(h,z))$. 
Its fundamental vector fields associated to the basis $\{u_{a+1}, \ldots , u_{c}\}$ are $\{ (X_{a+1}, W_{a+1}),
\ldots , (X_{c},W_c)\}$.
Given $h \in H$, we take $\Upsilon_h \colon K \times N  \to K \times 
N $ the map defined by $\Upsilon_h 
(k,z)  = \Upsilon(h,(k,z))$.
Then, we have
$$
\Ima 
\Pi^*  =
\left \{ \omega \in  {\bigwedge}^*  
(\gamma_{b+1}, \ldots,\gamma_f)
\otimes 
\lau{\Om}{*}{\per{p}}{N /  \mathcal{N}}\ \Bigg/ \
\begin{array}{l}
\hbox{(i) } i_{X_i} \omega = - i_{W_i}\omega \hbox{ if  } a < i \leq c \} \\[,3cm]
\hbox{(ii) }  L_{X_i} \omega =  -L_{W_i}\omega \hbox{ if  } a < i \leq c \},
\\[,3cm]
\hbox{(iii) } (\Upsilon_h)^*\omega = \omega \hbox{ for } h \in H.
\end{array}
\right\}.
$$
Let $H_{0}$  be the unity connected component of $H$. Recall that the 
subgroup $H_{0}$ is normal in $H$ and that the quotient 
$H/H_{0}$ is a finite group. Conditions (ii) gives that 
$\om$ is $H_{0}$-invariant. So,  condition (iii) can be replaced by:
$
\hbox{(iv) } (\Upsilon_h)^*\omega = \omega \hbox{ for } h \in H/H_{0}.
$
Therefore
$$
\Ima 
\Pi^*  =
\left \{ \omega \in  {\bigwedge}^*  
(\gamma_{b+1}, \ldots,\gamma_f)
\otimes 
\lau{\Om}{*}{\per{p}}{N /  \mathcal{N}}\ \Bigg/ \
\begin{array}{l}
\hbox{(i) } i_{X_i} \omega =  -i_{W_i}\omega \hbox{ if  } a < i \leq c \} \\[,3cm]
\hbox{(ii) }  L_{X_i} \omega = - L_{W_i}\omega \hbox{ if  } a < i \leq c \}.
\end{array}
\right\}^{H/H_{0}}.
$$
Since the group $H/H_{0}$ is a finite one, we get that the 
cohomology $\hiru{H}{*}{\Ima \Pi^* }$ is isomorphic to
$
\left(\hiru{H}{*}{A^\cdot}\right)^{H/H_{0}}, 
$
where $A^{*} $ is the differential complex
\bee
\left \{ \omega \in  {\bigwedge}^*  
(\gamma_{b+1}, \ldots,\gamma_f)
\otimes 
\lau{\Om}{*}{\per{p}}{N /  \mathcal{N}}\ \Bigg/ \
\begin{array}{l}
\hbox{(i) } i_{X_i} \omega =  -i_{W_i}\omega \hbox{ if  } a < i \leq c \} \\[,3cm]
\hbox{(ii) }  L_{X_i} \omega =  -L_{W_i}\omega \hbox{ if  } a < i \leq c \}
\end{array}
\right\}.
\eee
So, it remains to compute $\hiru{H}{*}{A^{\cdot}}$. This computation can be 
simplified by using these three facts:
\begin{itemize}
    \item[-] $i_{W_i}\omega =L_{W_i}\omega  =0$ for each $a <  i \leq b$,
    since the foliation $ \mathcal{N}$  is defined by the action of $G \cap H$.
    \item[-] $i_{X_i}\gamma_j = 
    \delta_{{ij}}$ for all
$ i ,j$ (cf.  \ref{tp}.\otto). 
\item[-] $d\gamma_j = 
    0$ for 
$
b< j$ (cf.  \eqref{ciclo}). 
\end{itemize}
We 
get that $A^*$ is the differential complex
$$
\left \{ \omega \in  {\bigwedge}^*  
(\gamma_{b+1}, \ldots,\gamma_f)
\otimes 
\lau{\Om}{*}{\per{p}}{N /  \mathcal{N}}\ \Bigg/ \
\begin{array}{l}
\hbox{(i) } i_{X_i} \omega = - i_{W_i}\omega \hbox{ if  } b < i \leq c \} \\[,3cm]
\hbox{(ii) }  0 =  L_{W_i}\omega \hbox{ if  } b < i \leq c \}
\end{array}
\right\}= 
$$
$$
{\bigwedge}^*  (\gamma_{c+1}, \ldots,\gamma_{f}) \otimes 
\underbrace{\left \{ \omega \in  {\bigwedge}^*  
(\gamma_{b+1}, \ldots,\gamma_c)
\otimes 
\lau{\Om}{*}{\per{p}}{N /  \mathcal{N}} \Bigg/  
\begin{array}{l}
\hbox{(i) } i_{X_i} \omega   =  - i_{W_i}\omega \hbox{ if  } b < i \leq c \} \\[,3cm]
\hbox{(ii) }  0  =   L_{W_i}\omega \hbox{ if  } b < i \leq c \}
\end{array}
\right\}}_{B^*}.
$$
 A straightforward computation gives that the canonical writing of a form $ \omega \in  {\bigwedge}^*  
(\gamma_{b+1}, \ldots,\gamma_c)
\otimes 
\lau{\Om}{*}{\per{p}}{N /  \mathcal{N}}$ verifying (i) is
\be
\label{writ}
\omega  = 
 \om_0 + \sum_{b < i_1 < \cdots < i_\ell\leq c} \ (-1)^\ell \gamma_{i_1} 
\wedge \dots 
\wedge \gamma_{i_\ell} \wedge 
(i_{W{_{i_\ell}} }\cdots i_{W{_{i_1}}}\om_0)
\ee
for some $\om_0 \in \lau{\Om}{*}{\per{p}}{N /  \mathcal{N}}$  (cf. Lemma \ref{iX}). 

Consider now $b<i,j\leq c$. Since $K/G$ is an abelian group (cf. Proposition \ref{trio}) and $H$ is a Lie group then ${\displaystyle [W_i,W_j] = \sum_{l=a+1}^b C^l_{ij} W_l}$. Then, $i_{[W_i,W_j]}\om_0 =0$ since 
the foliation $ \mathcal{N}$  is defined by the action of $G \cap H$. So, the canonical writing of a form $\om \in B^*$ is \eqref{writ}
for some $\om_0 \in  
 \left\{ \eta\in \lau{\Om}{*}{\per{p}}{N /  \mathcal{N}}\ /
\ L_{W_i} \eta= 0 \hbox{ if  } b < i \leq c \right\} = \left(\lau{\Om}{*}{\per{p}}{N / \mathcal{N}} \right)^{H_0} .
$

Then, the operator 
$
\Delta \colon B^* \TO \left(\lau{\Om}{*}{\per{p}}{N / \mathcal{N}} \right)^{H_0} ,
$
defined by
$
\Delta (\om) = \om_0$, is a 
 differential isomorphism.
 We
conclude that the differential complex $A^*$ is isomorphic to
$
{\bigwedge}^*  
(\gamma_{c+1}, \ldots ,\gamma_{f})
\otimes   \left(\lau{\Om}{*}{\per{p}}{N / \mathcal{N}} \right)^{H_0} 
$
 and therefore 
$
 \hiru{H}{*}{A^\cdot} 
\cong
\hiru{H}{*} {K / \mathcal{E} } 
\otimes
\lau{\IH}{*}{\per{p}}{N / \mathcal{N}}
$
(cf. \refp{KE} and \eqref{hache}). Since the operator $\Delta$ is $(H/H_0)$-equivariant (cf. \ref{tp}.\doce) then we get 
\begin{eqnarray*}
\lau{\IH}{*}{\per{p}}{K \bi{\times}{H} N / \mathcal{W} } &= &
\hiru{H}{*} {\Ima \Pi^*}=  \left(\hiru{H}{*} {A^\cdot}  \right)^{H/H_0} =
\left( \hiru{H}{*} {K / \mathcal{E} }  \otimes \lau{\IH}{*}{\per{p}}{N /  \mathcal{N}  }\right)^{H/H_{0}}.
\end{eqnarray*}
This ends the proof. \qed

\prg{\bf Remarks.}
        \addtocounter{Trece}{\value{proposition}} 

\Zati When the Lie group $G$ is commutative then  $K$ is also commutative. Differential forms $\gamma_{\small\bullet}$ are $K$-invariants on the left and on the right, so $\left(\hiru{H}{*} {K / \mathcal{E} }  \right)^{H} = \hiru{H}{*} {K / \mathcal{E} } $  and therefore 
\begin{equation*}
\lau{\IH}{*}{\per{p}}{K \bi{\times}{H} N / \mathcal{W} } = \hiru{H}{*} {K / \mathcal{E} }  \otimes \left( \lau{\IH}{*}{\per{p}}{N / \mathcal{N}}\right)^{H/H_{0}} =\hiru{H}{*} {K / \mathcal{E} } \otimes \left( \lau{\IH}{*}{\per{p}}{N / \mathcal{N}}\right)^{H}
\end{equation*}
as it has been proved in \cite[Proposition 3.8.4]{SW3}.

\zati  \addtocounter{Catorce}{\value{zatia}}
Since the foliation $\mathcal{E}$ is a riemannian foliation defined on a compact manifold then we know that the cohomology $\hiru{H}{*} {K / \mathcal{E} } $ is finite (cf. \cite{ESH}). So, the finiteness of $\lau{\IH}{*}{\per{p}}{K \bi{\times}{H} N / \mathcal{W} } $ depends on the finiteness of $ \lau{\IH}{*}{\per{p}}{N / \mathcal{N}}$.

\section{Finiteness of the BIC}
\label{fini}

We prove in this section that the BIC of a Killing foliation
 on a compact manifold
is finite 
dimensional. 
First of all, we present two geometrical tools we shall use in the proof: the isotropy type stratification and the Molino's blow up.

 We fix an isometric action $\Phi \colon G \times M \to M$ on the compact manifold $M$. We denote by
$\mathcal{F}$ the induced Killing foliation. For 
the study of $\mathcal{F}$ we can suppose that $G$ is connected (see 
Lemma \ref{normal}). We fix $K$ a  tamer group.  Notice that the group $G$ is normal in $K$ and the quotient $K/G$ is commutative (cf. Proposition \ref{trio}).

\prg{\bf Isotropy type stratification.} 
 \addtocounter{Cinque}{\value{proposition}} 
 The {\em isotropy type stratification} $\Str{K,M}  $ of $M$ is defined by the equivalence relation\footnote{For notions related with compact Lie group actions, we refer the reader to \cite{Br0}.}: 
\bee
x \sim y \Leftrightarrow K_x \hbox{ is conjugated to } K_y .
\eee
When $\depth \Str{K,M}  > 0$, any closed stratum $S \in \Str{K,M}  $ is a $K$-invariant submanifold of $M$ and then it possesses a $K$-invariant tubular neighborhood  $
(T,\tau,S,\R^m)$ whose structural group is $O(m)$. Recall that there are the following smooth maps associated with this
neighborhood:
\begin{itemize}
    \item[+] The {\em radius map} $\rho \colon T \to [0,1[$ defined fiberwise
    from the assignation
$[x,t] \mapsto t$. Each $t\not= 0$ is a regular value of the $\rho$.
The pre-image $\rho^{-1}(0)$ is $S$. This map is $K$-invariant, that is,
$\rho(k \cdot z) = \rho(z)$.
     \item[+] The {\em contraction} $H \colon T \times [0,1] \to T$
     defined fiberwisely from $([x,t],r ) \mapsto  [x,rt]$. The restriction
     $H_t \colon T \to T$ is an
     embedding for each
     $t\not= 0$   and $H_0 \equiv \tau$. We shall write $H(z,t) = t
     \cdot z$. This map is $K$-invariant, that is, $ t \cdot (k \cdot z) =
     k \cdot (t \cdot z)$.
\end{itemize}
The hyper-surface $D = \rho^{-1}(1/2)$ is the {\em
     tube} of the tubular neighborhood. It is a $K$-invariant submanifold of $T$.
     Notice that the map
$$
\nabla \colon D \times [0,1[ \TO T,
$$
defined
     by $\nabla(z,t) = (2t) \cdot z$ is a $K$-equivariant smooth map,
     where $K$ acts trivially on the $[0,1[$-factor. Its restriction
$\nabla \colon D \times ]0,1[ \TO T\menos S$ is a $K$-equivariant diffeomorphism.

 Denote $S_{_{min}}$ the union of closed (minimal) strata
and choose $
T_{_{min}}$ a disjoint family of $K$-invariant tubular neighborhoods  of the
 closed strata. The union of associated tubes is denoted by $D_{_{min}}$.
Notice that the induced map $\nabla_{_{min}} \colon D_{_{min}} \times ]0,1[ \TO
 T_{_{min}} \menos S_{_{min}}$ is a $K$-equivariant diffeomorphism.

\prg {\bf Molino's blow up}. 
  \addtocounter{Uno}{\value{proposition}} The Molino' blow up \cite{Mo} of the foliation $\mathcal{F}$ produces a
new foliation $\widehat{\mathcal{F}}$ of the same kind but of smaller depth.
 We suppose $\depth \Str{K,M}  > 0$.
 The {\em blow up} of $M$ is the compact manifold
$$
\widehat{M} =
\left\{
\Big( D_{_{min}} \times ]-1,1[\Big) \coprod \Big( (M\menos S_{_{min}}) \times \{
-1,1\}\Big)
\right\} \Big/ \sim ,
$$
where $(z,t) \sim (\nabla_{_{min}} (z, |t|), t/|t|)$, and the map
$
\mathcal{L} \colon  \widehat{M} \TO M
$
defined by
$$
\mathcal{L}(v) =
\left\{
\begin{array}{ll}
\nabla_{_{min}} (z, |t|)& \hbox{if } v = (z,t) \in  D_{_{min}} \times
]-1,1[ \\[,2cm]
z& \hbox{if }  v = (z,j) \in (M\menos S_{_{min}}) \times \{
-1,1\}.
\end{array}
\right.
$$
Notice that $\mathcal{L}$ is a continuous map whose restriction
$
\mathcal{L} \colon \widehat{M} \menos \mathcal{L}^{-1}(S_{_{min}}) \to M\menos S_{_{min}}
$
is a $K$-equivariant smooth trivial 2-covering.

Since the map $\nabla_{_{min}}$ is $K$-equivariant then $\Phi$ induces the action
$\widehat{\Phi} \colon K \times \widehat{M} \TO \widehat{M}$ by
saying that the blow-up $\mathcal{L}$ is $K$-equivariant.
The open
submanifolds $\mathcal{L}^{-1}(T_{_{min}})$  and
$\mathcal{L}^{-1}(T_{_{min}} \menos S_{_{min}})$
are clearly $K$-diffeomorphic to
$D_{_{min}} \times ]-1,1[$ and $D_{_{min}} \times (]-1,0[ \cup ]0,1[)$ respectively.

The
restriction
$\widehat{\Phi} \colon G \times \widehat{M} \TO \widehat{M}$
is an isometric action with $K$ as a tamer group.
The induced Killing foliation is $\widehat{\mathcal{F}}$.
Foliations $\mathcal{F}$ and $\widehat{\mathcal{F}}$ are
related by $\mathcal{L}$ which is a foliated map. Moreover, if $S$ is a not
minimal stratum of $\Str{K,M}  $ then there exists an unique  stratum $S'\in \Str{K,\wh{M}}$
such that $\mathcal{L}^{-1}(S) \subset S'$. The family $\{ S'  \ / \ S \in \Str{K,M}  \}$ covers $\widehat{M}$
and verifies the relationship: $S_1 \prec S_2 \Leftrightarrow S'_1 \prec
S'_2$. We conclude the important property
\be
\label{prof}
\depth \Str{K,\wh{M}} < \depth \Str{K,M}.
\ee

\prg{\bf Finiteness of a tubular neighborhood}.  We suppose $\depth \Str{K,M}  > 0$. Consider a closed stratum $S \in \Str{K,M}$.  Take $(T,\tau,S,\R^m)$ a $K$-invariant tubular neighborhood.  We fix a base point $x \in S$. The isotropy subgroup $ K_x$ acts orthogonally on the fiber $\R^m = \tau^{-1}(x)$. So, the induced action $\Lambda_x \colon G_x \times \R^m \to \R^m$ is an isometric action, it gives the Killing foliation $\mathcal{N}$ on $\R^m$.

\bP
\label{tub}
If the BIC of $(\R^m , \mathcal{N})$ is finite dimensional then the BIC of 
$(T,\mathcal{F})$ is also finite dimensional.
\eP

 \pro We proceed in two steps.
 
    \Zati \underline{\em $K_y=K_x$ \ for each $y \in S$}. 
     
 \smallskip
 
 The canonical projection $\pi \colon S \to S/K$ is an homogeneous bundle with fiber $K/K_x$. For any open subset $V \subset S/K$ the pull back  $\tau^{-1}\pi^{-1}(V)$ is a $K$-invariant subset of $T$, then we can apply the Mayer-Vietoris technics to this kind of subsets 
  (cf. \ref{bs}.\hspace{-.1cm} \quince).

 Since the manifold $S/K$ is a compact one then we can find a finite good covering $\{U_i \ / \ i \in I\}$ of it (cf. \cite{Br}). An inductive argument on the cardinality of $I$ reduces the proof of the Lemma  to the case where $T=\tau^{-1}\pi^{-1}(V)$, where $V$ is a contractible open subset of $S/K$.
 
 Here, the manifold $T$ is $K$-equivalently diffeomorphic to $V \times\left(K \bi{\times}{K_x} \R^m\right)$, where $K$ does not act on the first factor.
 So, the natural retraction of $V$ to a point gives a $K$-equivariant retraction of $T$ to the twisted product $K \bi{\times}{K_x} \R^m$. Now the result comes directly from   \ref{tp}.\hspace{-.1cm} \trece\catorce\ since  $(K,G,K_x)$ is a trio.

 \zati \underline{\em General case}. 
     
 \smallskip

The stratum $S$ is $K$-equivariantly diffeomorphic to the twisted product $K \bi{\times}{N(K_x)} F$ where $N(K_x)$ is the normalizer of $K_x$ on $K$ and $F = S^{K_x}$. So, the tubular neighborhood $T$  is $K$-equivariantly diffeomorphic to the twisted product $K \bi{\times}{N(H)} N$ where
$N$ is the  manifold  $\tau^{-1}(F)$. The previous case gives that the BIC of $(N,\mathcal{F}_N)$ is finite dimensional.
Now the result comes directly from   \ref{tp}.\trece\catorce\ since  $(K,G,N(K_x))$ is a trio.
\qed

The main result of this work is the following 
\bt
The BIC of the foliation determined by an isometric action on a compact 
manifold is finite dimensional.
\et

\pro
Let $\F$ be a Killing foliation defined on a compact manifold $M$ induced by an isometric action $\Phi \colon G \times M \to M$ where $G$ is a Lie group. Without loss of generality we can suppose that the Lie group $G$ is a connected one (cf. Lemma \ref{normal}). We fix a tamer group $K$. We know that $G$ is normal in $K$ and the quotient group $K/G$ is commutative (cf. Proposition \ref{trio}).
 
Let us consider the following statement
\begin{center}
    $\mathfrak{A}(U,\mathcal{F}) =$
``The BIC
    $
    \lau{\IH}{*}{\per{p}}{U/\mathcal{F}}
    $
    is finite dimensional for each perversity $\per{p}$,"
    \end{center}
    where $U \subset M$ is a $K$-invariant submanifold.
     We prove $\mathfrak{A}(M,\mathcal{F})$ by induction on $ \dim M $. The result is clear when $\dim M=0$. We suppose $\mathfrak{A}(W,\mathcal{F})$ for  any $K$-invariant compact submanifold $W$ of $M$ with $\dim W < \dim M$ and we prove $\mathfrak{A}(M,\mathcal{F}) $. We proceed in several steps.

\medskip

{\bf  First step: 0-depth.} Let us suppose $\depth \Str{K,M}=0$. Since $G \triangleleft K$ and  $K_x $ is conjugated to $ K_y$ then $G_x$ is conjugated to $ G_y$, $\forall x,y \in M$. We get that the foliation $\mathcal{F}$ is a
(regular) riemannian foliation (cf. \cite{Mo}). Its BIC is just the basic
cohomology (cf. \ref{bs}.\refbat). Then  $\mathfrak{A}(M,\mathcal{F})$ 
comes from \cite{ESH}.

\smallskip

{\bf  Second step: Inside \mbox{\boldmath $M$}.} Let us suppose $\depth \Str{K,M}>0$.
 The family $\left\{M \menos
 S_{_{min}},T_{_{min}}\right\}$ is a
basic covering of $M$ and the we get the exact sequence
(cf. \ref{bs}.\quince)
$$
0 \to
\lau{\Om}{*}{\per{p}}{\mf}
\to
\lau{\Om}{*}{\per{p}}{\left(M \menos
S_{_{min}}\right)/\mathcal{F}}
\oplus
\lau{\Om}{*}{\per{p}}{T_{_{min}}/\mathcal{F}}
\to
\lau{\Om}{*}{\per{p}}{\left(T_{_{min}} \menos
S_{_{min}}\right)/\mathcal{F}}
\to
0.
$$
The Five Lemma gives
\bee
\mathfrak{A}(T_{_{min}} \menos S_{_{min}},\mathcal{F}),  \ \
\mathfrak{A}(T _{_{min}},\mathcal{F}) \ \hbox{ and } \
\mathfrak{A}(M \menos S_{_{min}},\mathcal{F})
\Longrightarrow
\mathfrak{A}(M,\mathcal{F}).
\eee

Since   $T_{_{min}}\menos S_{_{min}}$ is $K$-diffeomorphic to $D_{_{min}} \times ]0,1[$ (cf. (cf. \ref{fini}.\cinque))
then 
$\mathfrak{A}(D_{_{min}},\mathcal{F})
\Longrightarrow
\mathfrak{A}(T_{_{min}}\menos S_{_{min}},\mathcal{F}).
$
The inequality $\dim D_{_{min}}< \dim M$ gives
\bee
\mathfrak{A}(T _{_{min}},\mathcal{F}) \ \hbox{ and } \
\mathfrak{A}(M \menos S_{_{min}},\mathcal{F})
\Longrightarrow
\mathfrak{A}(M,\mathcal{F}).
\eee

In order to prove  $\mathfrak{A}(T _{_{min}},\mathcal{F}) $ it suffices to prove $ \mathfrak{A}(T , \mathcal{F})$ where $(T,\tau,S,\R^m)$ a $K$-invariant tubular neighborhood of closed stratum $S$ of $\Str{K,M}$.
 Following
Proposition \ref{tub} we have
\bee
\mathfrak{A}(\R^m , \mathcal{N}) \Longrightarrow \mathfrak{A}(T , \mathcal{F}) \Longrightarrow \mathfrak{A}(T _{_{min}},\mathcal{F}).
\eee

Consider the orthogonal decomposition $\R^m = \R^{m_1} \times \R^{m_2}$, where $\R^{m_1} = \left( \R^m \right)^{G_x}$.
The only fixed point of the restriction $\Lambda_x \colon G_x \times \R^{m_2} \to \R^{m_2}$ is the origin. So, there exists a Killing foliation\footnote{It is given by the orthogonal action $\Lambda_x \colon G_x \times \S^{m_2-1} \to \S^{m_2-1}$.}  $\mathcal{G}$ on the sphere $\S^{m_2-1}$ with $(\R^{m_1} \times \R^{m_2} , \F) = (\R^{m_1} \times c\S^{m_2-1} ,\mathcal{I} \times c\mathcal{G})$. Propositions \ref{pro} and \ref{local5} give:
\bee
\mathfrak{A}(\S^{m_2-1} , \mathcal{G}) \Longrightarrow \mathfrak{A}(\R^{m_1} \times c\S^{m_2-1} ,\mathcal{I} \times c\mathcal{G})
\Longrightarrow \mathfrak{A}(\R^m , \mathcal{N}).
\eee
Finally, since $\dim \S^{m_2-1} < m \leq \dim T \leq \dim M$ we have
\be
\label{piso}
\mathfrak{A}(M \menos S_{_{min}},\mathcal{F})
\Longrightarrow
\mathfrak{A}(M,\mathcal{F}).
\ee

\smallskip

{\bf Third step: Blow-up.} Let us suppose $\depth \Str{K,M}>0$.
The family $\left\{\mathcal{L}^{-1}(M \menos
 S_{_{min}}),\mathcal{L}^{-1}(T_{_{min}})\right\}$  is a
basic covering of $\widehat{M}$ and the we get the exact sequence
(cf. \ref{bs}.\quince)
\bee
0 \to
\lau{\Om}{*}{\per{p}}{\widehat{M}/\widehat{\mathcal{F}}}
\to
\lau{\Om}{*}{\per{p}}{\mathcal{L}^{-1}(M \menos
S_{_{min}})/\widehat{\mathcal{F}}}
\oplus
\lau{\Om}{*}{\per{p}}{\mathcal{L}^{-1}(T_{_{min}})/\widehat{\mathcal{F}}}
\to
\lau{\Om}{*}{\per{p}}{\mathcal{L}^{-1}(T_{_{min}} \menos
S_{_{min}})/\widehat{\mathcal{F}}}
\to
0.
\eee
Following \ref{fini}.\uno\ we have that 
\begin{itemize}
\item[-] $\mathcal{L}^{-1}\left(M \menos
S_{_{min}}\right)$ is $K$-diffeomorphic to two copies of $M\menos S_{_{min}}$, 

\item[-]  $\mathcal{L}^{-1}\left(T_{_{min}} \right)$ is $K$-diffeomorphic to $D_{_{min}} \times ]-1,1[$,
\item[-]  $\mathcal{L}^{-1}\left(T_{_{min}} \menos S_{_{min}}\right)$ is $K$-diffeomorphic to $D_{_{min}} \times \left( ]-1,0[ \cup ]0,1[\right)$.

\end{itemize}
Now, the Five Lemma gives
$$
\mathfrak{A}(D_{_{min}},\widehat{\mathcal{F}} )\hbox{ and } \
\mathfrak{A}\left(\widehat{M},\widehat{\mathcal{F}}\right)
\Longrightarrow
\mathfrak{A}(M \menos
S_{_{min}},\mathcal{F}).
$$
But, the inequality $\dim D_{_{min}} < \dim M$ gives 
\be
\label{piso2}
\mathfrak{A}\left(\widehat{M},\widehat{\mathcal{F}}\right)
\Longrightarrow
\mathfrak{A}(M \menos
S_{_{min}},\mathcal{F}).
\ee

\smallskip

{\bf Forth step: Final blow-up.} When $\depth \Str{K,M}=0$ we get $\mathfrak{A}(M,\mathcal{F})$ from the First step. Let us suppose  $\depth \Str{K,M}>0$. From \eqref{piso} and \eqref{piso2} we get
\bee
\mathfrak{A}\left(\widehat{M},\widehat{\mathcal{F}}\right)
\Longrightarrow
\mathfrak{A}(M,\mathcal{F}).
\eee
with $\depth \Str{K,\widehat{M}} < \depth \Str{K,M}$ (cf. \eqref{prof}). By iterating this procedure we get
\bee
\mathfrak{A}\left(\wt{M},\wt{\mathcal{F}}\right) =\mathfrak{A}\left(\widehat{\mathop{\widehat{\mathcal{M}}}^{\cdots}},\widehat{\mathop{\widehat{\mathcal{F}}}^{\cdots}}\right) \Longrightarrow \cdots \Longrightarrow \mathfrak{A}\left(\widehat{M},\widehat{\mathcal{F}}\right)
\Longrightarrow
\mathfrak{A}(M,\mathcal{F}),
\eee
with $\depth \Str{K,\wt{M}} =0$. We finish the proof by applying again the First Step.\qed

\end{document}